\documentclass[12pt]{amsart}
\usepackage{bbm}
\usepackage{verbatim}
\usepackage{graphicx,color}
\usepackage{enumerate}
\usepackage[all]{xy}
\usepackage[hyperindex=true,plainpages=false,colorlinks=false,pdfpagelabels]{hyperref}

\theoremstyle{plain}
\newtheorem*{The*}{Theorem}

\newtheorem*{Cor*}{Corollary}

\theoremstyle{definition}

\newtheorem*{Rem*}{Remark}
\newtheorem*{Con*}{Conjecture}
\newtheorem{Con}{Conjecture}

\DeclareMathOperator{\Jac}{Jac}

\renewcommand{\Im}{\operatorname{Im}}
\renewcommand{\Re}{\operatorname{Re}}

\newcommand{\R}{\mathbb{R}}

\newcommand{\C}{\mathbb{C}}

\newcommand{\Z}{\mathbb{Z}}
\renewcommand{\P}{\mathbb{P}}

\begin{document}

\title{Towards a constrained Willmore conjecture}

\author{Lynn Heller}

\address{ Institut f\"ur Differentialgeometrie\\  Leibniz Universit{\"a}t Hannover\\ Welfengarten 1\\ 30167 Hannover\\ Germany
 }
 \email{lynn.heller@math.uni-hannover.de}
\author{Franz Pedit}

\address{ Department of Mathematics and Statistics\\
University of Massachusetts\\
Amherst, MA 01003\\
 USA}
 \email{pedit@math.umass.edu}



\date{\today}

\begin{abstract}
We give an overview of the constrained Willmore problem and address some conjectures arising from partial results and numerical experiments. Ramifications of these conjectures would lead to a deeper understanding of the Willmore functional over conformal immersions from compact surfaces.
\end{abstract}
\maketitle

\section{Introduction}
When bending an elastic membrane its bending energy per area can be shown to be  locally proportional to the square of its (mean) curvature. Very likely this has been known to Hooke as a manifestation of his {\em stress proportional to  strain law}. Mathematically we model the membrane by an immersed surface $f\colon M\to\R^3$. The bending energy density is then given by  $H^2 dA$ with $H$ the mean curvature  (the arithmetic mean of the principal curvatures) and $dA$ the area form both calculated with respect to the induced Riemannian metric 
$|df|^2$ on $M$. Thus, the total {\em bending energy} is
\[
\mathcal W(f)=\int_M H^2 dA
\]
which is also known as the {\em Willmore energy} of the immersion $f$. Indeed, in the mid 1960s Willmore \cite{Willmore}  proposed to study the minima and also the stationary points,  {\em Willmore surfaces}, of the functional $\mathcal W$ over immersions $f$ of compact surfaces $M$. He sought for an extrinsic analogue to the notion of tight immersions which are the minima of the total absolute Gaussian curvature $\int_M |K|dA$. Willmore showed that $\mathcal W(f)\geq 4\pi$ with equality attained at a round sphere of any radius since $\mathcal W(f)$ is invariant under scaling.  Furthermore, Willmore verified  that among tori of revolution in $\R^3$  the Clifford torus (and its scaled versions),  whose profile circle of radius $1$ has its center at distance $\sqrt{2}$ from the rotation axis,
is  the minimum with energy $2\pi^2$. He then posed the question, later to become the {\em Willmore conjecture}, whether 
the infimum of $\mathcal W$ over immersed surfaces of genus one is $2\pi^2$ and whether it is attained at the Clifford torus. 

Seemingly unknown to Willmore, in the 1920s Blaschke and his school \cite{Blaschke} studied surface theory in M\"obius geometry. Among other things they showed that one can associate to a surface $f\colon M\to \R^3$ in a M\"obius invariant way a {\em conformal Gauss map}  by attaching to any point $f(p)$ the sphere tangent to the surface of radius $1/H(p)$. Working in the light cone model of the conformal compactification $\R^3\cup\{\infty\}\cong S^3$, the space of oriented 2-spheres is parameterized by the Lorentzian 4-sphere $S^4_1$, a pseudo-Riemannian rank 1 symmetric space of constant curvature one. The M\"obius group of $3$-space is presented by the causality preserving isometries ${\bf O}^{+}(4,1)$ of $S^4_1$.  The M\"obius invariant induced area form of  the conformal Gauss map $\gamma_f\colon M\to S^4_1$ can be written as $(H^2-K)dA$. This motivated Blaschke and his school to study ``conformally area minimizing'' surfaces $f\colon M\to \R^3$,  that is the critical points of the functional  
\[
\tilde{\mathcal{W}}(f) =\int_M (H^2-K)dA
\]
We note that on compact surfaces of genus $g$  the functionals $\mathcal W$ and $\tilde{\mathcal W}$ differ by the Euler number $\int_M K dA= 4\pi(1-g)$ and thus have the same critical points: Willmore surfaces are Blaschke's conformal minimal surfaces and therefore invariant under the full group of M\"obius transformations of 3-space. 

Some of the, mostly local, results obtained by Blaschke and his collaborators have eventually been rediscovered or reproven:
\begin{itemize}
\item
The Euler-Lagrange equation of $\mathcal{W}$ is given by  $\triangle H+2H(H^2-K)=0$. Thus M\"obius transforms of minimal surfaces  are examples of Willmore surfaces;
\item
A surface $f\colon M\to \R^3$ is Willmore if and only if its conformal Gauss map $\gamma_f\colon M\to S^4_1$ is (branched) minimal;
\item 
Isothermic Willmore surfaces are (locally) M\"obius congruent to a minimal surface in one of the space forms \cite{Thomsen}.
\end{itemize}
The integrant $H^2-K=|\mathring{S}|^2$ is just the squared length of the tracefree 2nd fundamental form of the immersion $f\colon M\to \R^3$, which is known to scale inverse proportionally to the induced area  under M\"obius transformations (in fact, under any conformal changes of the ambient metric). Hence $\tilde{\mathcal W}=\int_M |\mathring{S}|^2 dA$ is defined for any conformal target manifold. In particular, ignoring the Euler number, we have 
$\mathcal W=\int_M(H^2 \pm 1)dA$ for immersions into the sphere $S^3$ and hyperbolic space $H^3$.

In the early 1980s Bryant \cite{Bryant_Willmore}  classified all Willmore 2-spheres $f\colon S^2\to \R^3$ as M\"obius inversions of minimal surfaces with planar ends in $\R^3$. The Willmore energy $\mathcal W=4\pi k$ is quantized by the number of ends and, except for $k=2,3,5,7$, all values are realized. For $k=1$ one obtains the absolute minimum of $\mathcal W$, a round 2-sphere or, after an inversion, a plane. For higher values of $\mathcal W$  there are smooth families of Willmore spheres. 

The classification of Willmore tori and a resolution of the Willmore conjecture turned out to be more involved.  The first lower bound for $\mathcal W$ over immersions $f$  of a Riemann surface $M$ into the n-sphere  $S^n$ was given in early 1980 in a paper by Li and Yau \cite{LiYau}. Due to  M\"obius invariance of the Willmore energy  one has
\[
\mathcal{W}(f)=\int_M(H^2+1)dA \geq  \sup_T \mathcal{A}(T\circ f)
\]
Here $T$ ranges over all M\"obius transformations of $S^n$ and $\mathcal{A}$ denotes the area functional. Therefore, one can estimate the minimal Willmore energy for a fixed conformal class from below by the conformal area 
$\inf_{f} \sup_T \mathcal{A}(T\circ f)$, where $f$ ranges over all conformal immersions of $M$ into $S^n$.
Applying a first eigenvalue estimate for the Laplacian, Li and Yau could show that for a
torus $M=\R^2/\Z\oplus \Z\tau$ whose  conformal structure lies in the domain $|\tau| \geq 1$, $\Im \tau \leq 1$, $0\leq \Re \tau \leq 1/2$ the conformal volume is bounded from below by $2\pi^2$.  In their study of minimal immersions given by first eigenfunctions of the Laplacian, Montiel and Ros \cite{MontielRos} (and independently Bryant \cite{Bryant-conf}) enlarged this domain slightly. They also provided examples of immersions of tori of rectangular conformal types whose
 conformal areas were below $2\pi^2$, thus erasing hopes that this approach would resolve the Willmore conjecture.  
Li and Yau also showed the useful estimate that if an immersion has a multiple point of order $k$ then $\mathcal W(f)\geq 4\pi k$. 

The existence of a minimizer for $\mathcal W$ over immersions of tori was first proven by Simon \cite{Simon} in the early 1990s, and then in early 2000 by Bauer and Kuwert  \cite{BauerKuwert} for any genera and in any codimension. But whether the minimizer in genus one  had to satisfy $\mathcal W\geq 2\pi^2$ remained open. Besides the cases covered in \cite{LiYau} the Willmore conjecture had been verified for channel tori (envelopes of a circle worth of round 2-spheres) \cite{JerPin} and tori invariant under an antipodal symmetry \cite{Ros}. Numerical evidence for the validity of the conjecture was obtained for Hopf Willmore tori \cite{Pinkall}, which were also the first examples of Willmore surfaces not M\"obius congruent to minimal surfaces, and for $S^1$-equivariant Willmore tori \cite{FerusPedit}. 

In the early 1990s the emerging integrable system techniques used to study harmonic maps from tori into symmetric spaces  \cite{Hi}, \cite{PS}, \cite{FerPedPinSte_S4}, \cite{BurFerPedPin_Annals} began to be applied to the conformal Gauss maps of Willmore surfaces. 
Even though a deeper appreciation of the complexity of the problem was obtained, for instance in the unpublished work by Schmidt \cite{Schmidt} in early 2000,
and new lower bounds in terms of the integrable systems structure were found, notably the quaternionic Pl\"ucker formula \cite{BLFPP}, this approach has so far failed to provide a resolution of the Willmore conjecture. It was not until early 2012 when Marques and Neves \cite{MarNev} gave a proof of the Willmore conjecture using very different  techniques based on a refinement of the Almgren and Pitts Min-Max approach. 

Knowing the existence of minimizers for any genus and the minimal value for $\mathcal W$ over immersions $f\colon M\to \R^3$ of tori, we can ask what to expect in higher genus. So far the only known examples of higher genus $g\geq 2$ Willmore surfaces in $\R^3$  are stereographic projections of compact minimal surfaces in the 3-sphere. The most prominent examples are Lawson's minimal surfaces $\xi_{g,1}$  of any genus $g$ with dihedral symmetries \cite{Lawson}, which naturally lend themselves as candidates for the minimizers \cite{Kusner}. Their Willmore energies $\mathcal W$ are increasing with the genus, starting from the Clifford torus at $\mathcal W=2\pi^2$ and limiting to $8\pi$ as the genus tends to infinity. 
There is some experimental evidence \cite{Kusneretal}, using bending energy decreasing flows, that Lawson's surfaces are indeed the minima for $\mathcal W$ on higher genus surfaces.

In addition to the existence results of minimizers in arbitrary codimension \cite{Simon, BauerKuwert} there has been renewed interest in the study of Willmore surfaces in higher codimension. In the late 1980s Ejiri \cite{Ejiri}  gave a description of Willmore 2-spheres admitting a dual Willmore 2-sphere in any codimension in terms of holomorphic lifts to a suitable twistor space. In codimension two this was made explicit in a paper by Montiel \cite{Montiel}, who showed that every Willmore 2-sphere in codimension two is either the inversion of a minimal surface with planar ends in $\R^4$ or the projection of a rational curve under the Penrose twister fibration $\C\P^3\to S^4$. The Willmore energy is quantized by $\mathcal{W}=4\pi k$ and all $k\in\Z$ occur. This result holds verbatim for Willmore tori in codimension two, provided the normal bundle has non-zero degree \cite{LeschkePeditPinkall}. Applying loop group techniques Dorfmeister and Wang \cite{DorfWang} recover some of those results in their language and also provide a general setting for studying Willmore surfaces in the framework introduced in \cite{DPW}. A more geometric approach, in the spirit of Calabi's characterization \cite{Calabi} of minimal 2-spheres, is given by Ma and his collaborators \cite{Ma}, \cite{MaWang}, \cite{MaWangWang}. This approach studies generalized Darboux transforms of Willmore surfaces and aims to construct a sequence of transforms which terminate, at least for Willmore 2-spheres, in a minimal surface with planar ends or a twistor projection of a rational curve.  Even though the existence of minimizers is known for any genus and in any codimension, it is still an open question whether the Willmore conjecture holds in higher codimension. It is conceivable, though hard to imagine,  that a minimizing torus in $\R^4$  has Willmore energy below $2\pi^2$. 

A refinement of the Willmore problem arises if we fix the conformal structure of the domain $M$ and look for minimizers (and more generally critical points) of the Willmore functional $\mathcal W$ under variations preserving the conformal type of $M$. The critical points are called (conformally)
{\em constrained Willmore surfaces}.  From the perspective of the integrable systems approach working in a fixed conformal class is more natural, if not required, since the conformal type is one of the integrals of motion.  There is also some hope that by utilizing Riemann surface theory - especially the description of conformally immersed surfaces in 3- and 4-space by a Dirac equation with potential \cite{BLPP} -  the 4th order non-linear elliptic analysis of the Euler-Lagrange equation can, to some extend, be replaced  to some extend by  the 1st order linear analysis of the Dirac equation. Such ideas may be found in the unpublished work of Schmidt \cite{Schmidt}, \cite{SchmidtCW} but have yet to be made precise.  

The conformal constraint introduces, as a Lagrange multiplier, a quadratic holomorphic differential $Q\in H^0(K^2)$ into the Euler-Lagrange equation
\begin{equation*}
\triangle H+2H(H^2-K)= \,<Q, \mathring{S}>
\end{equation*}
where $\mathring{S}$ denotes the tracefree 2nd fundamental form. Choosing $Q=0$ shows that every Willmore surface is a constrained Willmore surface.
 For a constant mean curvature surface we choose  $Q$ to be its holomorphic  Hopf differential to  
 verify that those surfaces (and their M\"obius transforms) are constrained Willmore. Since there are no holomorphic differentials on a genus zero Riemann surface constrained Willmore spheres are Willmore spheres.  There is a subtlety when deriving the constrained Euler-Lagrange equation: the subspace of conformal immersions of a Riemann surface inside the smooth manifold of all immersions is singular exactly at isothermic immersions. Since many of the known examples are isothermic (for instance, constant mean curvature tori) the validity of the equation at those surfaces required a deeper analysis of second variations  \cite{SchätzleIsotherm}.

 Existence of constrained minimizers in any conformal class has been proven for any genus and any codimension \cite{KuwSch_const}, \cite{KuwertLi}, \cite{Riv} provided that the given conformal class can be realized by an immersion of Willmore energy below $8\pi$. This condition is used to avoid that the minimizer acquires branch points (a branch point - as a limit of a double point - would imply that $W\geq 8\pi$ by the estimate of Li and Yau). In particular, the constrained minimizers with  Willmore energy below $8\pi$ are embedded tori.  In addition, the Willmore energy $\mathcal W$ varies continuously over the conformal types of those minimizers \cite{KuwSch_const}. Numerical experiments with bending energy decreasing flows \cite{HPSW}, \cite{CPS} suggest that constrained minimizers of genus one, whose conformal structures are sufficiently far from rectangular structures, will have $W\geq 8\pi$.  If one allows codimension at least two, then every conformal type \cite{LammSchätzle} can be realized by a conformal immersion with Willmore energy $\mathcal W \leq 8\pi$.  There is some evidence \cite{SchmidtCW} that in codimension one the upper  bound is $12\pi$. 
 
On the other hand, independent of the proof of the Willmore conjecture \cite{MarNev}, Ndiaye and Sch\"atzle \cite{NdiayeSchätzle1}
show that for rectangular conformal types $\R^2/\Z\oplus ib\Z$ sufficiently near to the square structure 
the minimizers in $\R^3$  are the (stereographic projection of) homogeneous flat tori $S^1\times S^1(1/b)$ in 3-spheres. They also generalize this result to arbitrary codimension  \cite{NdiayeSchätzle2}.  

After this somewhat broad stroked and kaleidoscopical overview of the subject, we will discuss some of the emerging conjectures concerning constrained Willmore surfaces together with supporting evidence and partial results. 

{\em Acknowledgments\,}: part of this research was supported by  the German Science Foundation DFG. The first author is indebted to the Baden-W\"urttemberg Foundation for supporting her research through the ``Elite Program for Postdocs". The second author was additionally supported by a University of Massachusetts  professional development grant. All images were produced by Dr. Nicholas Schmitt using his Xlab software suite. 

\section{The constrained Willmore conjecture}

\begin{figure}\label{fig:torus-tree}
\centering
\includegraphics[width=0.5\textwidth]{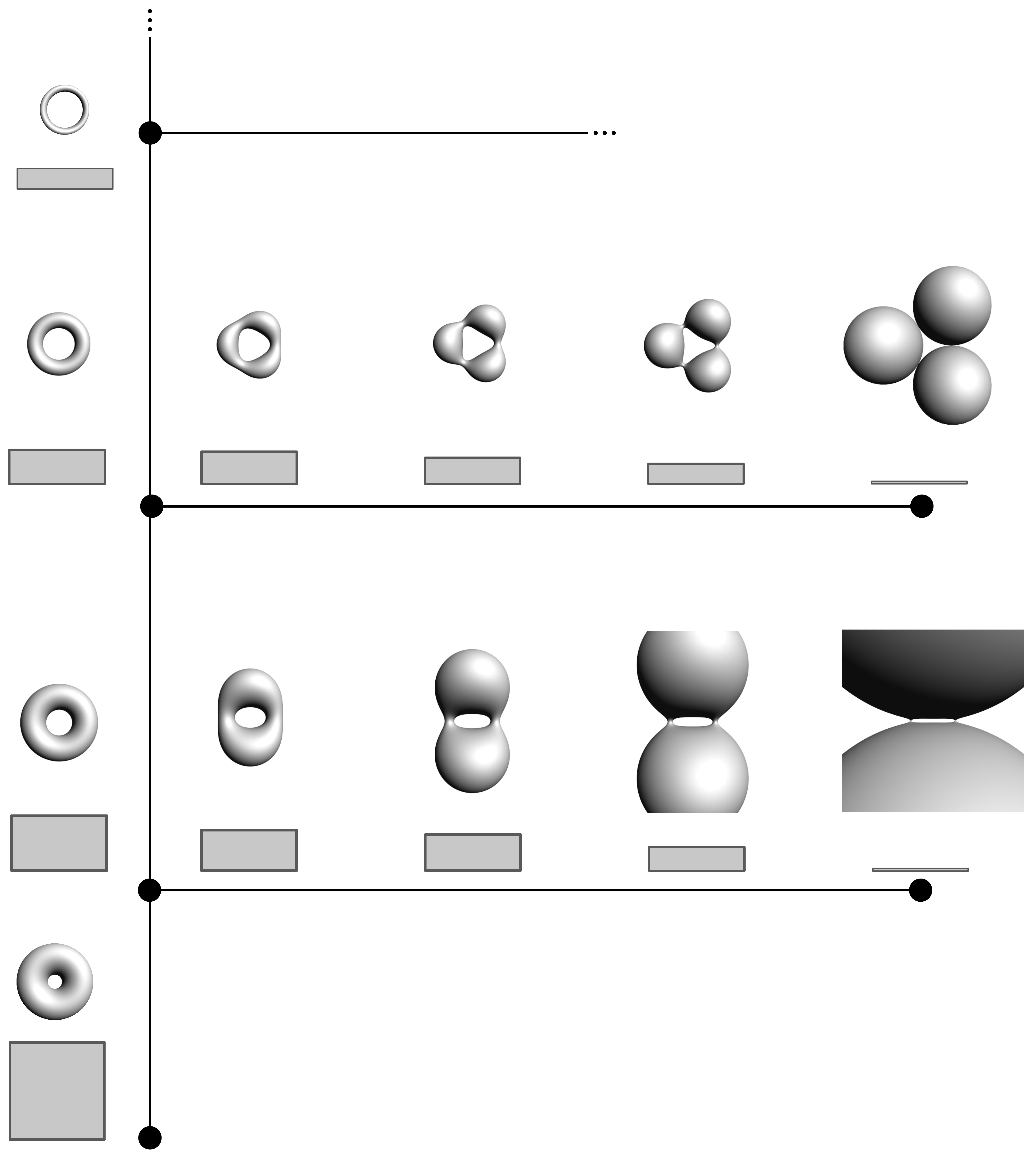}
\caption{
The vertical stalk represents the family of homogenous tori,
starting with the Clifford torus at the bottom.
Along this stalk are bifurcation points at which
the embedded Delaunay tori appear along the horizontal lines.
The rectangles indicate the conformal types.
}
\end{figure}

The classical Willmore conjecture asserts that among immersions of tori in 3-space the Clifford torus is the absolute minimum with Willmore energy $2\pi^2$. The Clifford torus  is also the unique embedded minimal torus in the 3-sphere, a result proven recently by Brendle \cite{Brendle} and conjectured in the 1970s by Lawson
\cite{Lawson}.  Since surfaces of constant mean curvature are constrained Willmore, it might be  reasonable to guess that the constrained minimizers in a given conformal class could be found among constant mean curvature tori.  We already saw that the existence of constrained minimizers \cite{KuwSch_const}  is guaranteed only for surfaces with $\mathcal W <8\pi$ (in which case they have to be embedded). Hence one should first consider embedded constant mean curvature tori as candidates.  By a recent result of Andrews and Li \cite{AndrewsLi} based on ideas in \cite{Brendle}, such surfaces have to be rotational tori in the 3-sphere. This result had been conjectured by Pinkall and Sterling \cite{PS} as a generalization of Lawson's conjecture. 
Embeddedness restricts those surfaces to be the products of circles, the homogeneous tori $S^1\times S^1(b)$, for any $b\geq 1$, and the unduloidal $k$-lobed Delaunay tori \cite{KilianSchmidtSchmitt1} in a 3-sphere for $k\geq 2$ (see Figure \ref{fig:torus-tree}). The conformal types $\R^2/\Z\oplus\tau \Z$ of those tori 
sweep out all rectangular conformal types $\tau=ib$  for $b\in[1,\infty)$ starting at the square structure.  The $k$-lobed Delaunay tori branch off from homogeneous tori at the conformal type  $\tau=i\sqrt{k^2-1}$ and limit to a bouquet of $k$ spheres of degenerate conformal type. The Willmore energy 
$\mathcal{W}(S^1\times S^1(b))=\pi^2 (b+1/b)$ of homogenous tori is unbounded and hence they cannot minimize $\mathcal W$ for all rectangular conformal types. 
In fact, Kuwert and Lorenz \cite{KuwLor} calculate  the Jacobi operator for $\mathcal{W}$ along homogeneous tori and show that negative eigenvalues appear precisely at $b=\sqrt{k^2-1}$ for $k\geq 2$, where the Delaunay tori branch off. Combining arguments in \cite{KilianSchmidtSchmitt1}  and \cite{KilianSchmidtSchmitt2} the Willmore energy  along the family of homogeneous tori $S^1\times S^1(b)$ for  $1\leq b\leq \sqrt{k^2-1}$ and continuing along the unduloidal $k$-lobed Delaunay tori is monoton increasing between $2\pi^2\leq \mathcal W < 4\pi k$ (the latter  being the Willmore energy of the bouquet of $k$-spheres). 

Thus, considering the case $k=2$, we have an embedded constant  mean curvature torus with $2\pi^2\leq  \mathcal W < 8\pi$ for every rectangular conformal type $\R^2/ \Z\oplus ib\Z$ which minimizes $\mathcal{W}$ among constant mean curvature tori in those conformal classes. Taking into account the results of  Ndiaye and Sch\"atzle \cite{NdiayeSchätzle1}, \cite{NdiayeSchätzle2}, who showed that for rectangular conformal types near the square structure the homogeneous tori $S^1\times S^1(b)$ minimize, we are led to a constrained Willmore conjecture for rectangular conformal types: 

\begin{Con}[2-lobe conjecture]
\label{conj:2-lobe}
The constrained minimizers of the Willmore energy  for tori in $\R^3$  of rectangular conformal types $\R^2/\Z\oplus ib\Z$ are (stereographic projections of) the homogeneous tori $S^1\times S^1(b)$ in the 3-sphere  for $1\leq b \leq \sqrt{3}$ and the 2-lobed Delaunay tori in a 3-sphere  for $b > \sqrt{3}$ limiting to a twice covered equatorial 2-sphere as $b\to \infty$ (see Figure \ref{fig:torus-tree}). 
\end{Con}

\begin{figure}\label{fig:candidates}
\vspace{0.5cm}
\includegraphics[width= 0.3\textwidth]{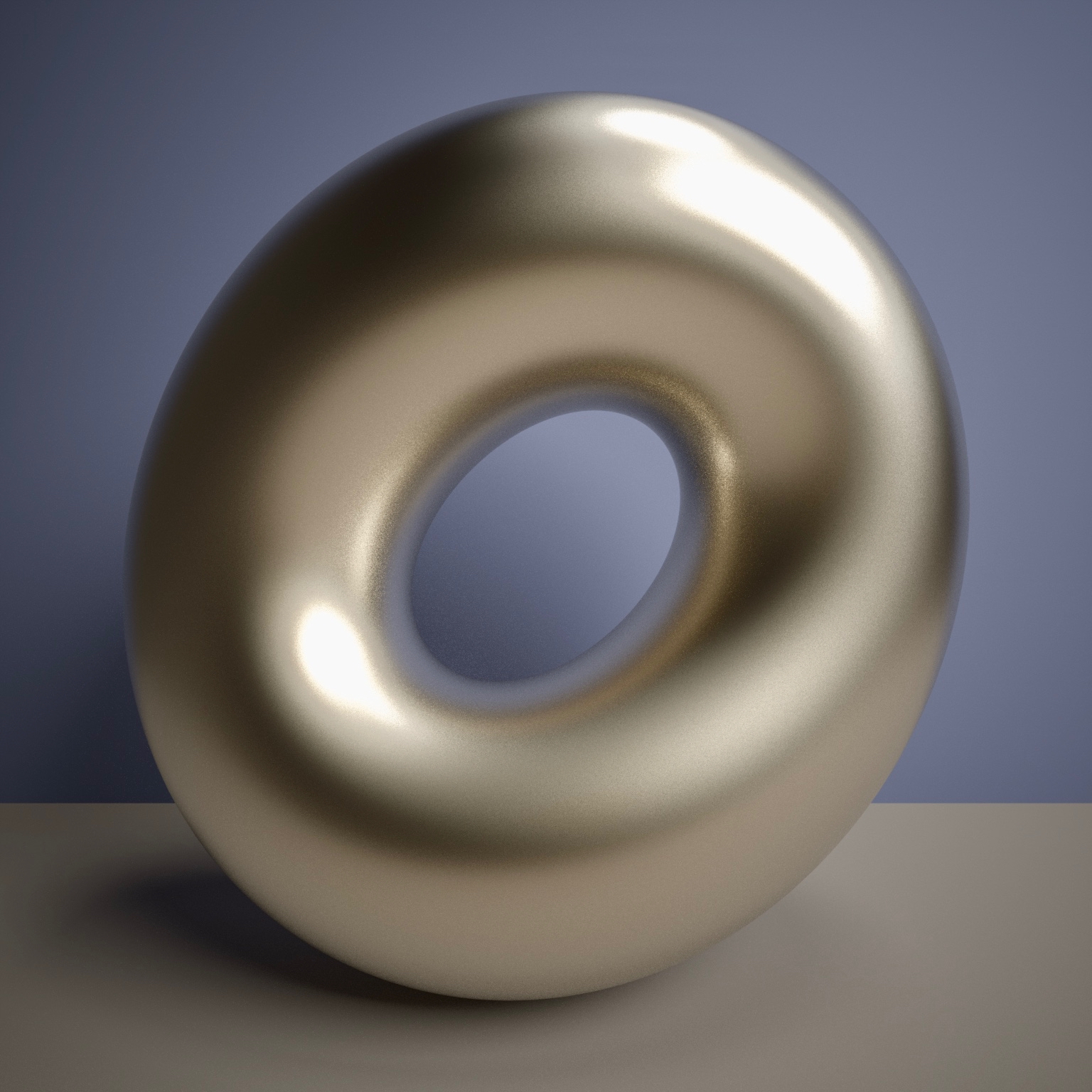}
\includegraphics[width= 0.3\textwidth]{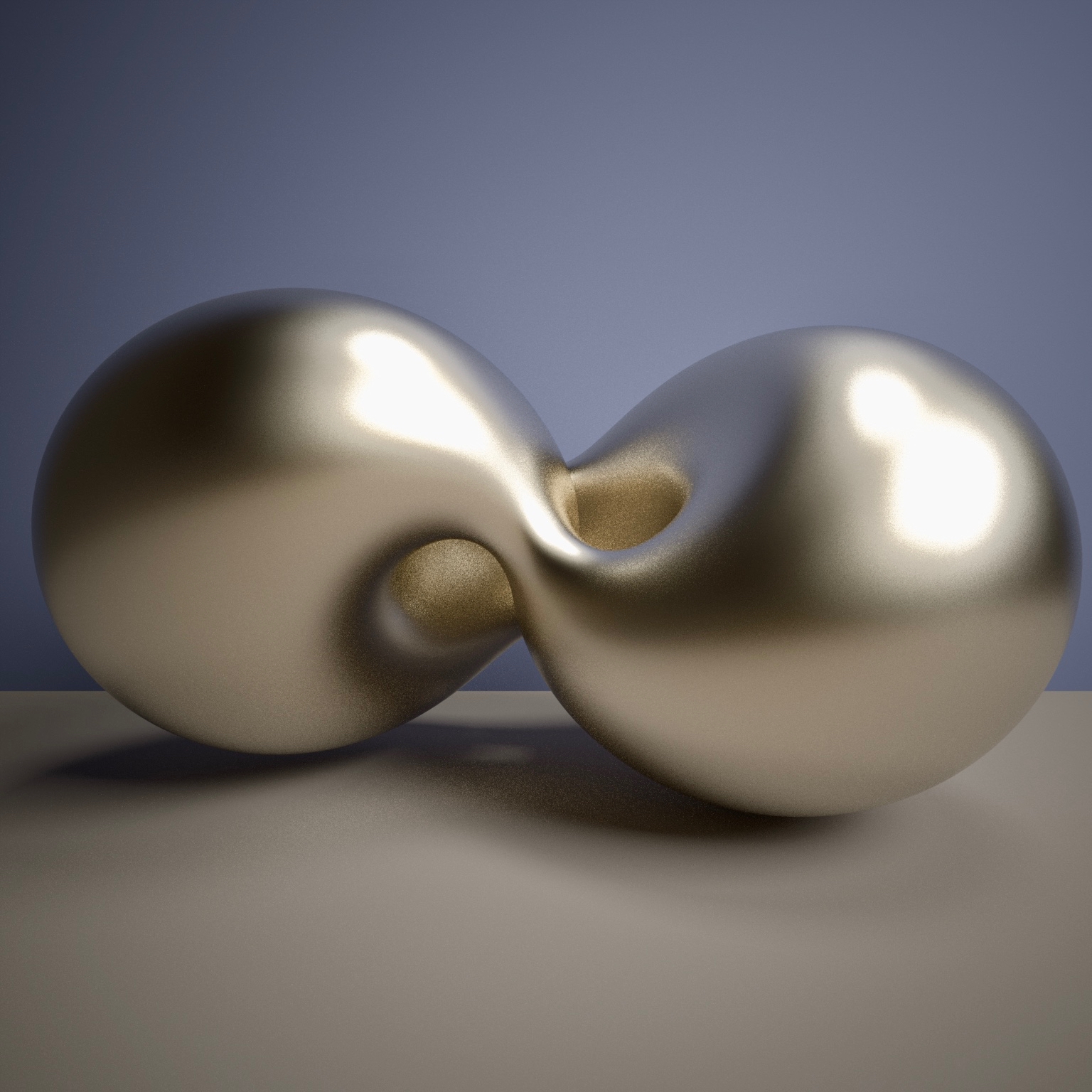}\\
\vspace{0.1cm}
\includegraphics[width= 0.3\textwidth]{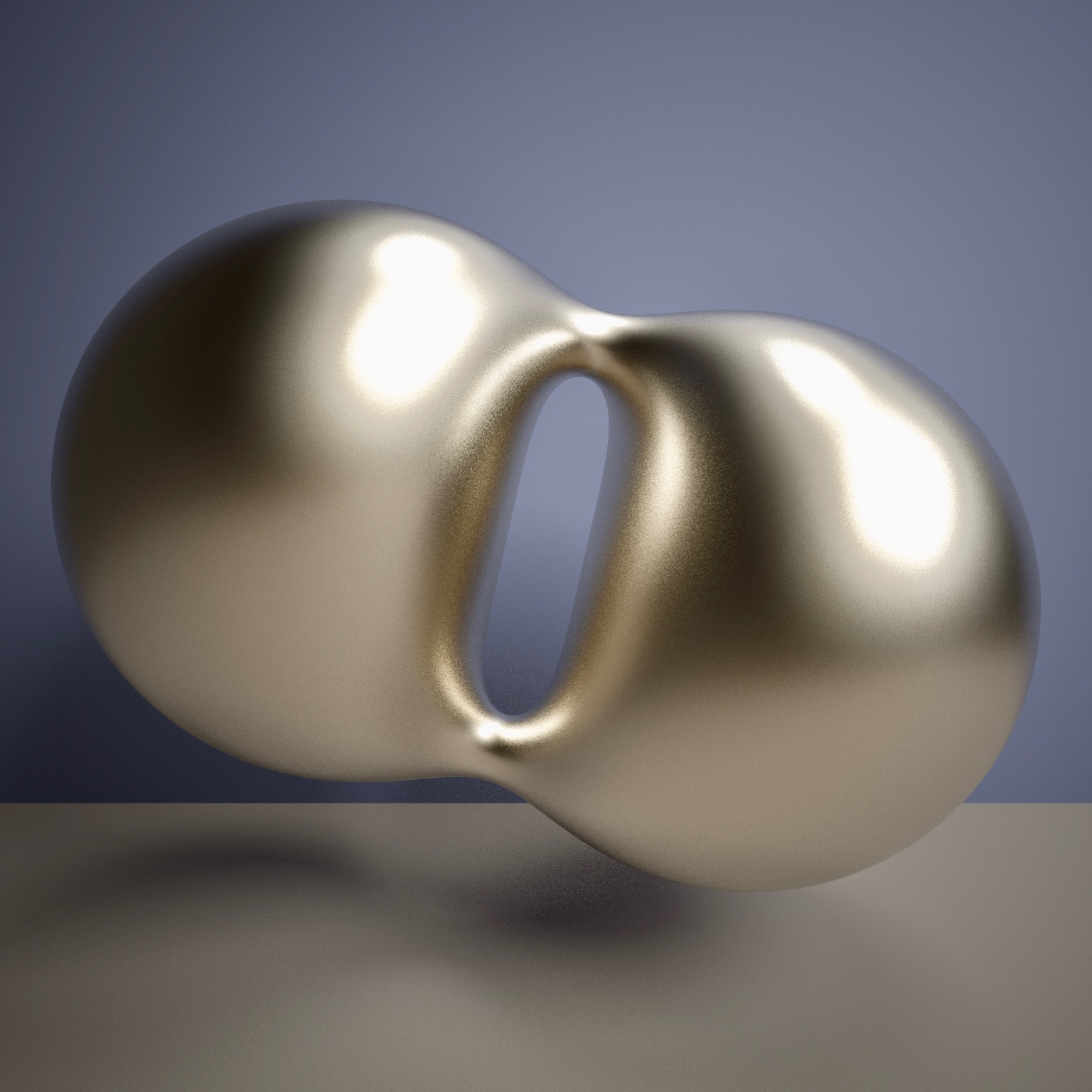}
\includegraphics[width= 0.3\textwidth]{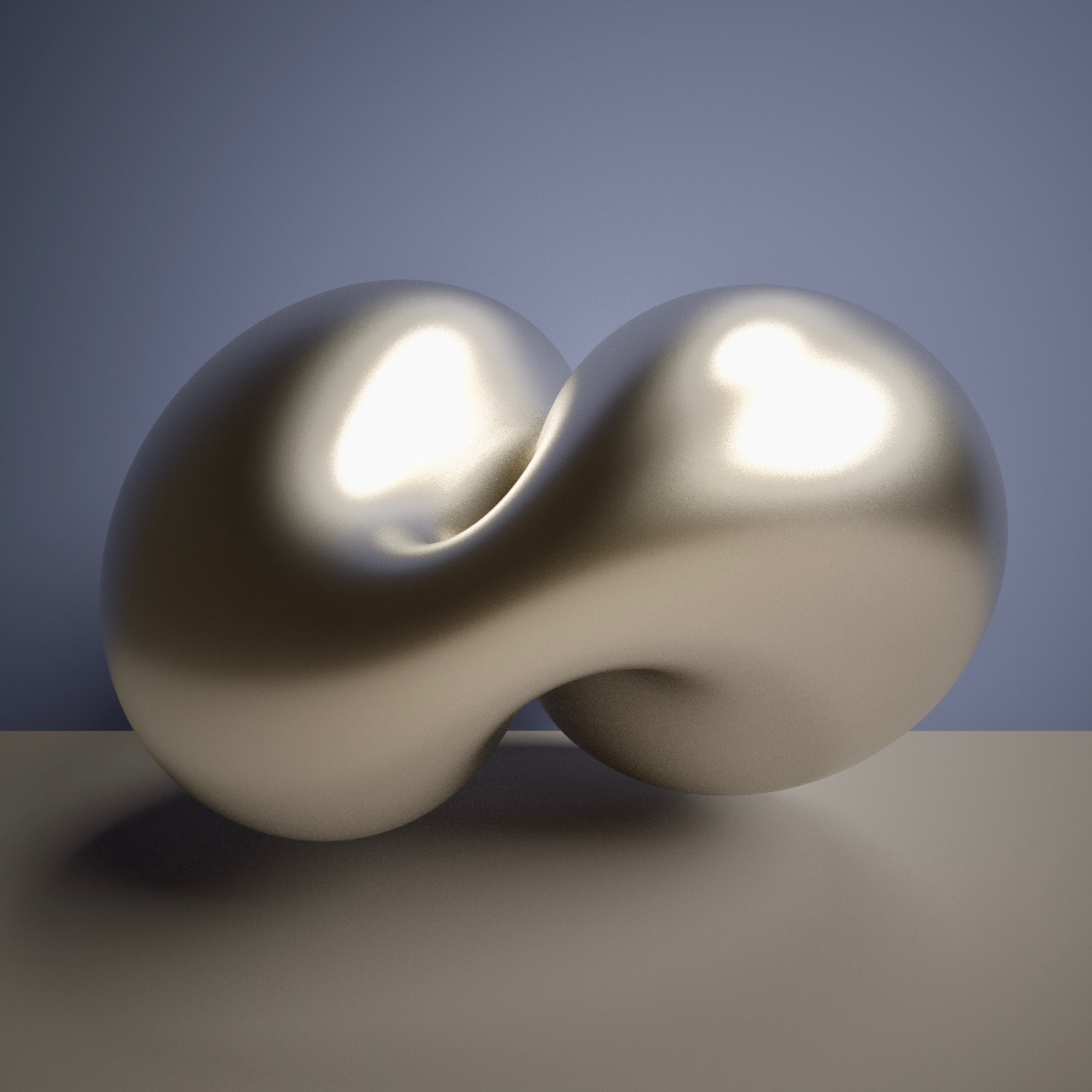}
\caption{ 
Experimentally found minimizers of the constrained Willmore problem \cite{HPSW}. They resemble $(1,2)$-equivariant constrained Willmore tori with intrinsic period $1$ - compare to Figure  \ref{fig:1lobe}. The conformal structures $\tau=a+ib$ of the shown tori are as follows:  upper row has $b < \sqrt{3}$, lower row has 
$b > \sqrt{3}$; left column has $a\sim 0$, right column has $a>>0$.  The torus on the bottom right is not embedded.} 
\end{figure}

\begin{figure}\label{fig:1lobe}
\vspace{0.5cm}
\includegraphics[width= 0.3\textwidth]{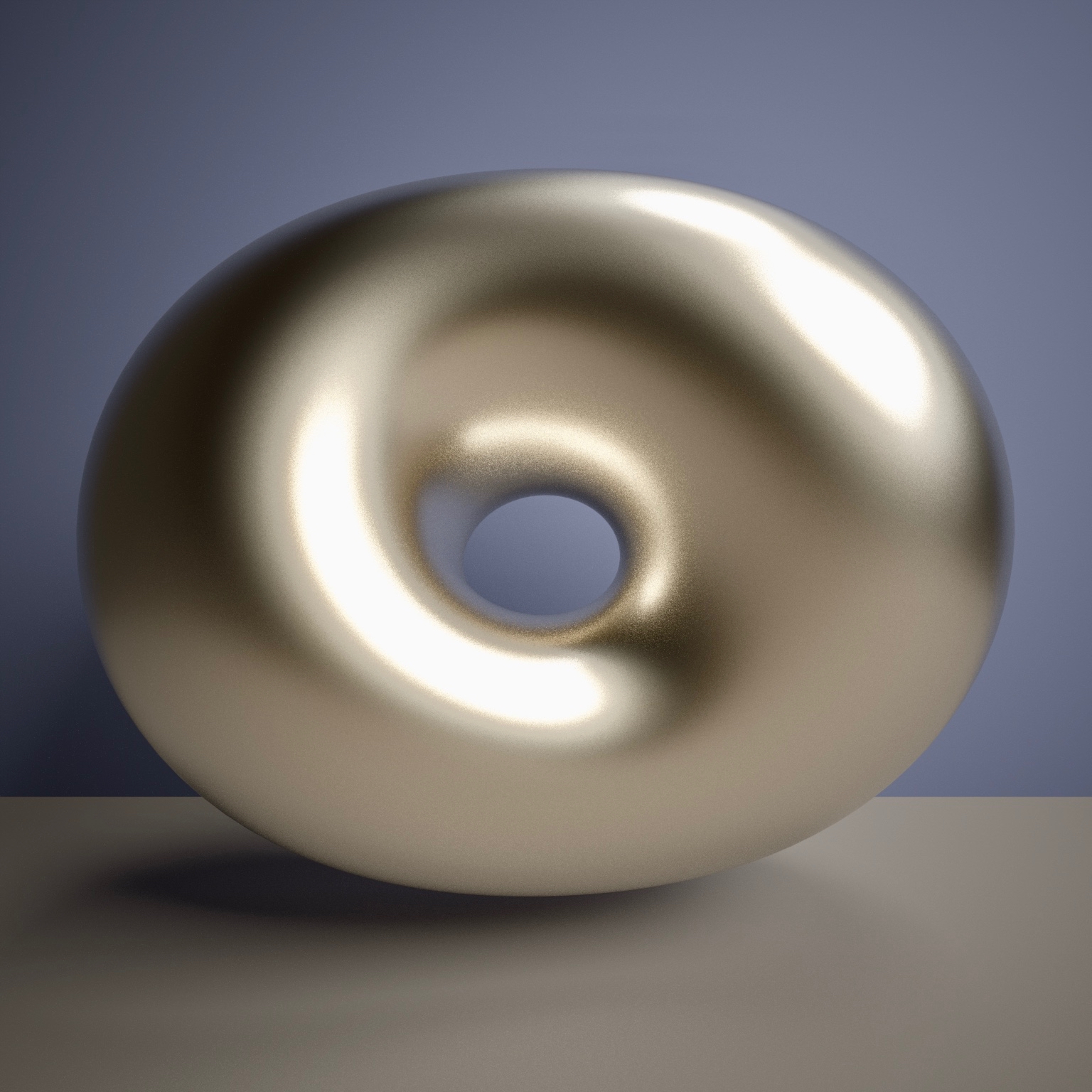}\hspace{1cm}
\includegraphics[width= 0.3\textwidth]{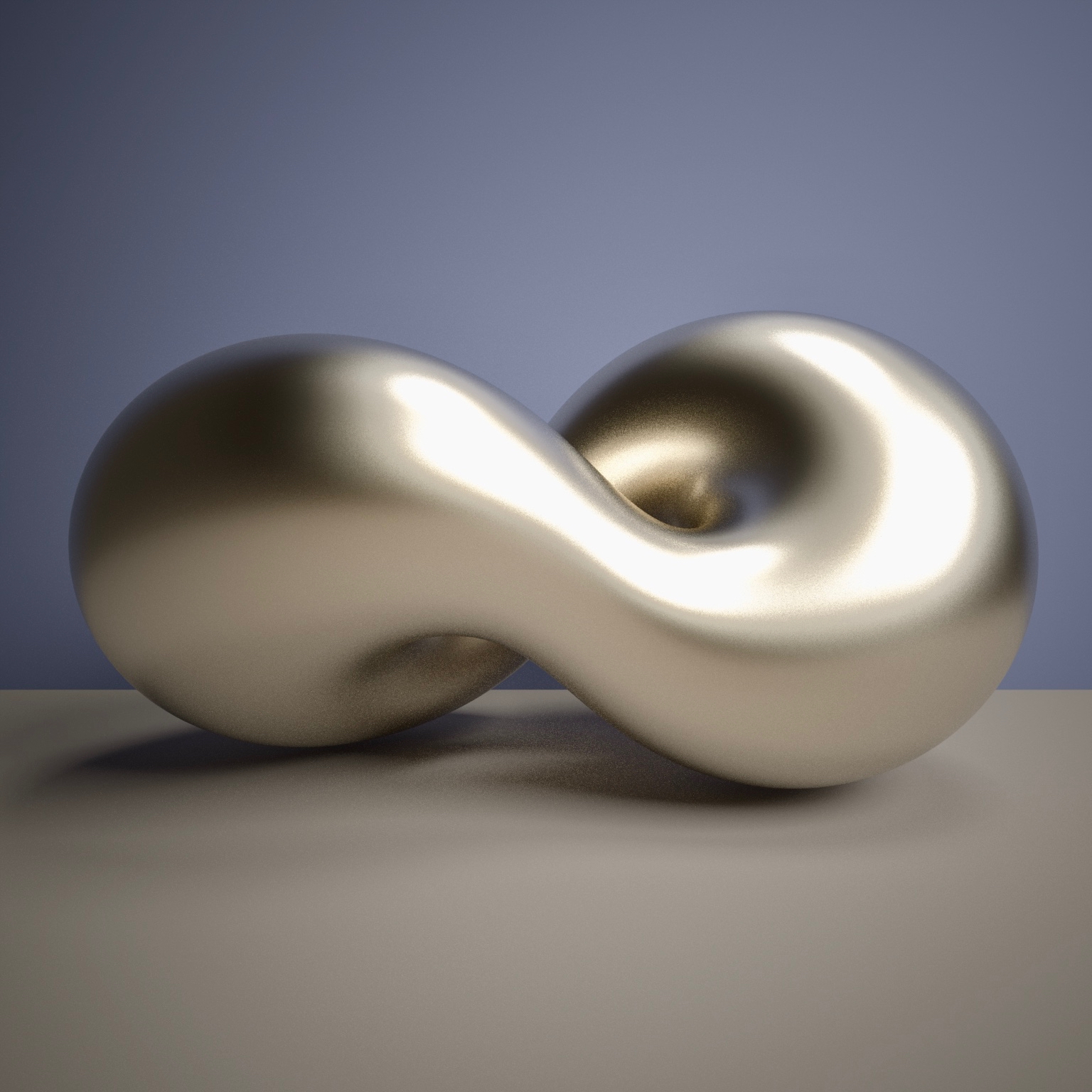}
\caption{ 
$(1,2)$-equivariant constrained Willmore tori with intrinsic period $1$. The tori lie in a $2$ parameter family of constrained Willmore tori deforming the Clifford torus and minimize the Willmore functional in conformal classes near the square structure.}
\end{figure}

For constrained minimizers in non-rectangular conformal classes there is much less guidance due to a lack of examples. 
Constrained minimizers, whose non-rectangular conformal types lie in a sufficiently small open neighborhood of the square structure will, by continuity, still have Willmore energy below $8\pi$ and thus have to be embedded. Therefore, they cannot have constant mean curvature by our previous discussion. The first examples of constrained Willmore tori {\em not} M\"obius congruent to constant mean curvature tori were constructed by Heller \cite{Hel1}, \cite {Hel2},  who described all 
 $S^1$ equivariant constrained Willmore tori in the 3-sphere where the orbits of the action are $(m,n)$ torus knots. 
 For instance, $(1,1)$-equivariance yields constrained Hopf Willmore tori which arise as preimages under the Hopf projection $S^3\to S^2$  of length and enclosed area constrained elastica on the round 2-sphere of curvature $4$. Their conformal types $\tau=a+ib$ are determined by the enclosed areas $a$ and lengths $b$ of the elastica and all conformal types occur. In a recent preprint 
Heller and Ndiaye \cite{NdiHel} generalize the results in  \cite{NdiayeSchätzle1} and characterize the minima of the constrained Willmore problem in a certain neighborhood of rectangular conformal classes: 
\begin{The*}\label{thm:LynnNdiaye}
For every $b \sim 1$ and $ b \neq 1$ there exists a sufficiently small $a(b)>0$ such that for every $a \in [0,a(b))$  the $(1,2)$ equivariant constrained Willmore torus $f_{a,b}$  of intrinsic period $1$ (see Figure~\ref{fig:1lobe}) and conformal type  $\tau=a+ib$ is a constrained Willmore minimizer. Moreover, for fixed $b \neq 1$ and  $b \sim 1$ the energy profile $a\mapsto \mathcal{W}(f_{a,b})$ is concave and varies real analytically over  $[0, a(b))$.
\end{The*}
It is expected  that this result can be extended to a neighborhood of all rectangular conformal types. In fact, laying a grid over the space of conformal structures and using an implementation of the conformal Willmore flow \cite{CPS}, we saw evidence that the $(1,2)$-equivariant constrained Willmore tori minimize \cite{HPSW} in all conformal classes (see Figure \ref{fig:1lobe}). These experiments also suggest that for conformal structures sufficiently far from rectangular types these tori have Willmore energy  $\mathcal{W}\geq 8\pi$ and in general will not be embedded (see Figure \ref{fig:candidates}).
 We emphasize though that for such conformal structures there is, as of yet, no existence proof guaranteeing immersed minimizers.  
 Based on these observations we are led to  

\begin{Con}[Constrained Willmore Conjecture]
The minimizers of the Willmore functional $\mathcal W$ over tori in a prescribed conformal class are  given by the  $(1,2)$-equivariant constrained Willmore tori (of intrinsic period $1$).  
\end{Con}
In contrast to the results in \cite{NdiayeSchätzle2} for rectangular conformal classes near the square structure,  the above theorem cannot be generalized to higher codimension:  it relies on the stability computations in \cite{KuwLor}  for codimension $1$. We do expect to find other minimizers in higher codimension: for instance the minimal torus of hexagonal conformal type in the 5-sphere \cite{MontielRos} of Willmore energy $\mathcal{W}= 4\pi^2/\sqrt{3}< 8\pi$. In fact, as already mentioned,  in codimension at least $2$ every conformal type \cite{LammSchätzle} can be realized by a conformal immersion with Willmore energy $\mathcal W \leq 8\pi$.  

\section{The constrained Willmore Lawson conjecture}
In the above discussions we have seen that having a rich enough class of examples at one's disposal helps guiding the theory.  Understanding $(1,2)$-equivariant tori and comparing them to the experimentally found minimizers gave some credence to the constrained Willmore conjecture. 
Unfortunately, so far the only known examples of compact constrained Willmore surfaces are equivariant,  or surfaces M\"obius congruent to constant mean curvature surfaces. To obtain a more detailed understanding of constrained Willmore surfaces a wider class of examples seems desirable. Since the conformal Gauss map $\gamma_f\colon M\to S^{4}_{1}$   of a  Willmore surface $f\colon M\to \R^3$  is a harmonic  (in fact, branched minimal) map into the space of round 2-spheres, one can apply the integrable systems theory  for harmonic maps into symmetric target spaces to construct  examples of  Willmore tori. 

The associated family of the harmonic map $\gamma_f$  into the Lorentz 4-sphere $S^{4}_{1}$ gives rise to a $\C^{\times}$ family of flat ${\bf sl}(4,\C)\cong {\bf so}(5,\C)$ connections \cite{DPW}, \cite{DorfWang}, \cite{FerPedPinSte_S4}
\[
\nabla^{\lambda}=\nabla+\lambda^{-1}\Phi + \lambda \bar{\Phi}
\]

Here $\Phi\in \Omega^{1,0}\otimes{\bf sl}(4,\C)$ is given by the $(1,0)$ part $\Phi=d\gamma_f^{1,0}$ of the derivative of $\gamma_f$ and $\nabla$ is the pullback of the Levi-Civita connection on $S^4_1$. There is a reality condition, namely for $\lambda\in S^1\subset \C^{\times}$ the $\nabla^{\lambda}$ have to be ${\bf sp}(1,1)$ connections. In this reformulation of the Euler-Lagrange equation for Willmore surfaces one can augment the loop $\nabla^{\lambda}$ by 
a quadratic holomorphic differential $Q$ and obtain a verbatim characterization of constrained Willmore surfaces as a family of flat connections 
\cite{BurQui}, \cite{Boh_Will}. The conformal Gauss map $\gamma_f$  then is a constrained branched minimal map.  In the case when  $M=T^2$ is a 2-torus the flat connection $\nabla^{\lambda}$ has four distinct parallel line subbundles except at isolated values of $\lambda\in\C^{\times}$ where pairs of the line bundles coalesce. Bohle \cite{Boh_Will} has shown that there can only be finitely many such points and that the line bundles have limits as $\lambda\to 0,\infty$. Thus, there exists a compact Riemann surface  $\lambda\colon \Sigma\to \P^1$, quadruply covering $\P^1$, parametrizing the parallel line bundles. The Riemann surface  $\Sigma$ is called the {\em spectral curve} and its genus $g_s$ the {\em spectral genus} of the constrained Willmore torus.  Now a point $L\in\Sigma$ presents a parallel line bundle for $\nabla^{\lambda(L)}$ and the family of connections $\nabla^{\lambda}$ is holomorphic in $\lambda$. Therefore, for fixed $p\in T^2$,  the fibers  $L_p$, as $L$ varies over $\Sigma$,  fit to a holomorphic line bundle $L(p)\to\Sigma$.  The complete integrability of constrained Willmore tori is encoded in the statement that the map 
\[
T^2 \to \Jac(\Sigma)\colon p\to L(p)
\]
is a group homomorphism - linear flow - of prescribed direction from the 2-torus $T^2$ into the (real) Jacobi torus $\Jac(\Sigma)$ of holomorphic line bundles over $\Sigma$. Since a group homomorphisms of a given slope  is determined by its value at one point $L(p_0)$, we can associate to a Willmore torus $f\colon T^2\to \R^3$ its {\em  spectral data}: the spectral curve $\Sigma$ and an initial condition $L(p_0)\in \Jac(\Sigma)$  for the linear flow. It has been shown in \cite{BLPP}, \cite{BPP} that one can reconstruct the constrained Willmore torus up to M\"obius equivalence from its spectral data. Notice that once $\Sigma$ is chosen there is a $g_s=\dim \Jac(\Sigma)$ dimensional freedom of choosing the initial condition $L(p_0)$ of the flow. This accounts for  {\em isospectral} deformations of a constrained Willmore torus: two of those dimensions conformally reparametrize the surface, but the $g_s-2$ dimensions transverse to the linear flow $T^2\to \Jac(\Sigma)$ yield non M\"obius congruent constrained Willmore tori (necessarily of the same conformal type and Willmore energy).  

What has been described here is a geometric manifestation of the finite gap theory of  the Novikov-Veselov hierarchy: the isospectral deformations account for the higher flows of the hierarchy. In the special case of constant mean curvature tori the spectral curve 
$\Sigma\to \P^1$ is hyperelliptic.
This case, corresponding to the $\sinh$ (or $\cosh$) Gordon hierarchy, has been studied in great detail \cite{PS}, \cite{Hi}, \cite{Bob} and many examples of higher spectral genus $g_s\geq 2$ constant mean curvature tori are known (for instance,  the Wente and Dobriner tori). 

The spectral genus $g_s$ gives some information about the complexity of the corresponding constrained Willmore torus \cite{Hel1}, \cite{Hel2}: 
\begin{itemize}
\item 
$g_s=0$ characterizes the homogeneous tori $S^1\times S^1(b)$ for $b\in[1,\infty)$; 
\item
$g_s=1$ characterizes equivariant constant mean curvature tori constructed in \cite{KilianSchmidtSchmitt1}; 
\item
 $g_s=2$ characterizes surfaces M\"obius congruent to  constant mean curvature tori of Wente type in any of the three space forms, and equivariant constrained Willmore tori in the associated family of constrained Hopf Willmore cylinders. In particular, the conjectured constrained minimizers in a given conformal class, the $(1,2)$-equivariant constrained Willmore tori, are found here. 
\end{itemize}
In general, equivariant constrained Willmore tori of orbit type $(m,n)$ have spectral genus $g_s\leq 3$ and there is a complete classification \cite{Hel1}, \cite{Hel2} of these surfaces for the case $g_s\leq 2$, including explicit parametrizations.
 
The lowest spectral genus where one can find constrained Willmore tori, which are neither equivariant nor M\"obius congruent to constant mean curvature tori, has to be at least $g_s\geq 3$, but its value is unknown. 
The reconstruction of a constrained Willmore torus form spectral data is in principle possible \cite{BLPP}, \cite{BPP} even though the  technical details have not been worked out. To construct such new examples of constrained Willmore tori  and to understand their structure would significantly increase our confidence in 
\begin{Con}[Constrained Willmore Lawson conjecture]
An embedded  constrained Willmore torus is equivariant.
\end{Con}
We know that the  conjecture holds for an embedded constrained Willmore torus which is M\"obius congruent to a constant mean curvature torus by the results of Brendle \cite{Brendle} and Li and Andrews \cite{AndrewsLi}.  An affirmative answer to this conjecture has significant ramifications:
\begin{itemize}
\item
{\em  The classical Willmore conjecture}: since the (unconstrained) minimizer of the Willmore functional has to be embedded, it would have to be equivariant. But the Willmore conjecture holds (or at least can be verified) for equivariant Willmore tori.
\item 
{\em  The 2-lobe conjecture}:  since the constrained minimizers in rectangular conformal classes have to be embedded, it would suffice to study the equivariant case to resolve the constrained Willmore conjecture for rectangular conformal types. 
\end{itemize}

\section{The stability conjecture}
We are now entering into much less charted terrain. There is some evidence that Willmore energy  decreasing flows on immersed 2-spheres limit to the round sphere. This has been shown using long time existence of the Willmore gradient flow  on 2-spheres with $\mathcal W\leq 8\pi$ in \cite{KuwSch}, \cite{KuwSch_removable} and is suggested by experiments in \cite{CPS} using a conformal Willmore flow. Therefore one is led to believe that the only stable Willmore spheres are round 2-spheres.  In what follows stability always refers to the positive (semi) definitness $\delta^2\mathcal{W}\geq 0$  of the Jacobi operator for the Willmore functional.  

Although there is  less supporting theoretical evidence for surfaces of genus one, we formulate 
\begin{Con}[Stability conjecture]
The only stable Willmore torus is the Clifford torus. 
\end{Con}
As we have explained in the previous section, every Willmore torus $f\colon M\to \R^3$ has some spectral genus $g_s$ and is contained in the $g_s-2$ dimensional isospectral family of non M\"obius congruent Willmore tori. If $V$ denotes the kernel (modulo M\"obius tranformations) of the Jacobi operator $\delta^2 \mathcal W_f$, then we have
\[
\dim V \geq g_s - 2
\]
Therefore,  a {\em strictly stable} Willmore torus $f\colon M\to \R^3$, that is $f$ is stable and  $V$ is trivial, has spectral genus $g_s \leq 2$.
All Willmore tori of spectral genus $g_s \leq 2$ are known \cite{Hel2} to be either M\"obius congruent to a minimal torus of Wente type in one of the space forms, or associated to the Hopf Willmore surfaces found by Pinkall \cite{Pinkall}. These surfaces could in principle be checked case by case to affirm the stability conjecture for strictly stable Willmore tori. 

There is further support for the conjecture involving the multiplicity  $\mu_1$ (modulo  M\"obius transformations) of the first eigenvalue $\lambda_1$ of the Jacobi operator $\delta^2\mathcal W_f$.  Since stability implies $\lambda_1\geq 0$, we have either $\dim V=0$ (for $\lambda_1 > 0$) or $\dim V = \mu_1$ (for $\lambda_1=0$). Therefore a stable Willmore torus has spectral genus at most $g_s\leq \mu_1+2$. For example, if one could show (as is the case for the Laplacian) that $\mu_1=1$ then the stability conjecture could be verified by checking only those Willmore tori of spectral genus $g_s=3$. 

Among Hopf Willmore tori  \cite{Pinkall} the stability conjecture follows from the fact that the analogous conjecture holds for elastic curves, the critical points of the total squared curvature:  stable elastica on the 2-sphere have to be great circles \cite{LanSin}. Thus, the only stable Hopf Willmore torus \cite{Pinkall} is the Clifford torus.

The following weaker version of the stability conjecture, which is of interest in its own right, would already provide an alternative proof of the Willmore conjecture:
\begin{Con}[Weak stability conjecture]
A stable Willmore torus has to be isothermic.
\end{Con}
Assuming the conjecture to hold, Thomsen's  result \cite{Thomsen} implies that an isothermic Willmore torus is (locally) M\"obius congruent to a minimal surface in one of the space forms. If the space form is $S^3$ then, due to the fact that a Willmore minimizer is embedded, we obtain an embedded minimal torus. In this case Brendle's resolution \cite{Brendle} of the Lawson conjecture \cite{Lawson} implies that the minimizer has to be the Clifford torus.
The other possibility is that the Willmore torus is an inversion of a minimal torus with planar ends in $\R^3$. But those surfaces are never embedded \cite{BT}. Finally, the Willmore torus could be one of the minimal tori coming from doubly periodic solutions of the $\cosh$ Gordon equation: these tori consist of two minimal pieces contained in the upper respectively lower hyperbolic half spaces $H^3$ smoothly connecting across the boundary at infinity, along which the torus has umbilic curves. The known examples \cite{BaBob} are not embedded but it is unknown whether all such tori arising from $\cosh$ Gordon solutions have self intersections. 

\begin{figure}\label{fig:highergenus}
\vspace{0.5cm}
\includegraphics[width= 0.4\textwidth]{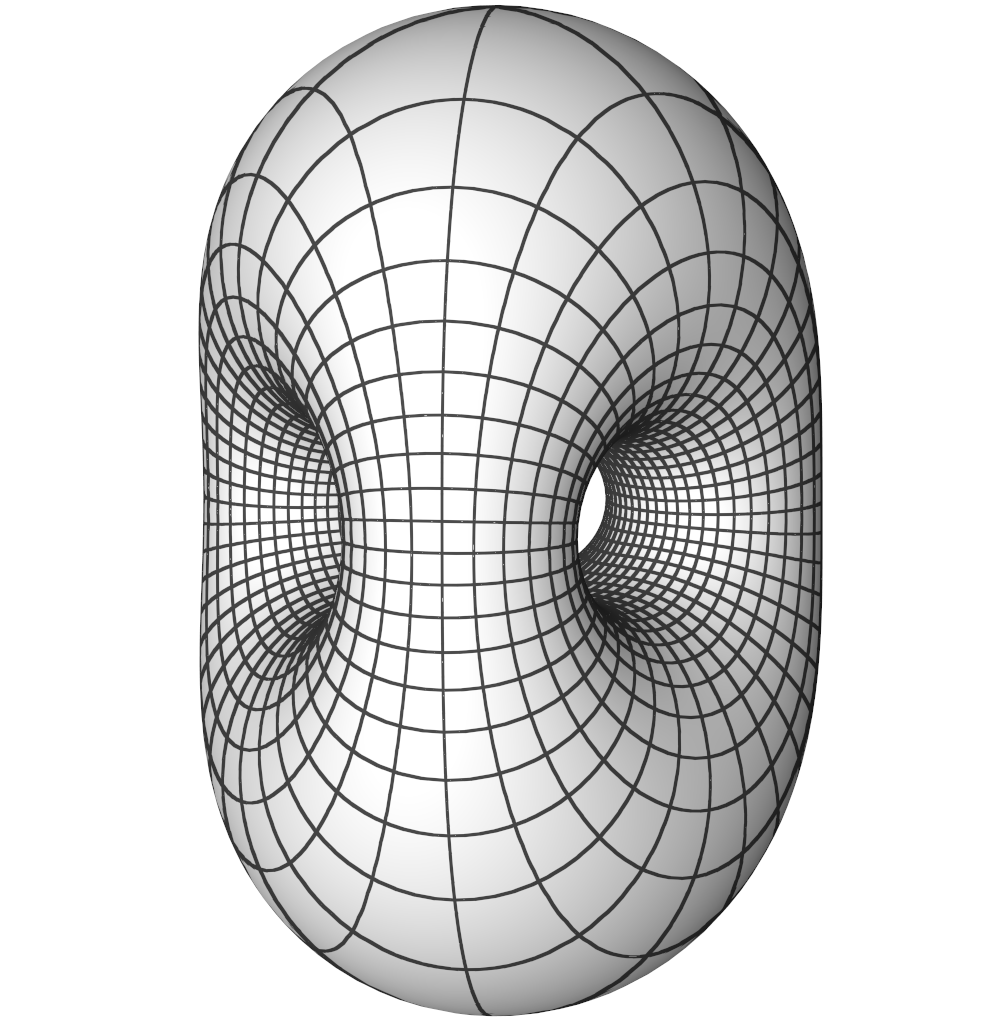}
\includegraphics[width= 0.4\textwidth]{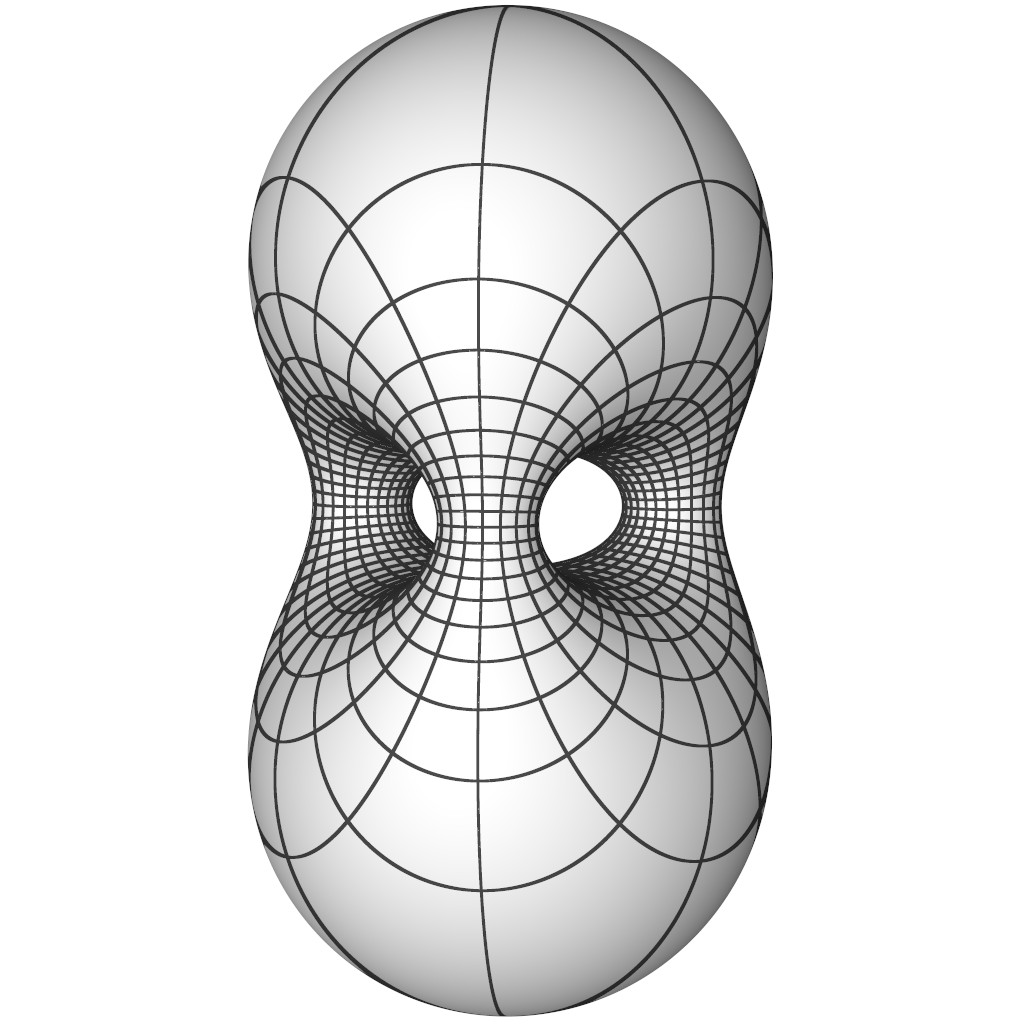}
\includegraphics[width= 0.4\textwidth]{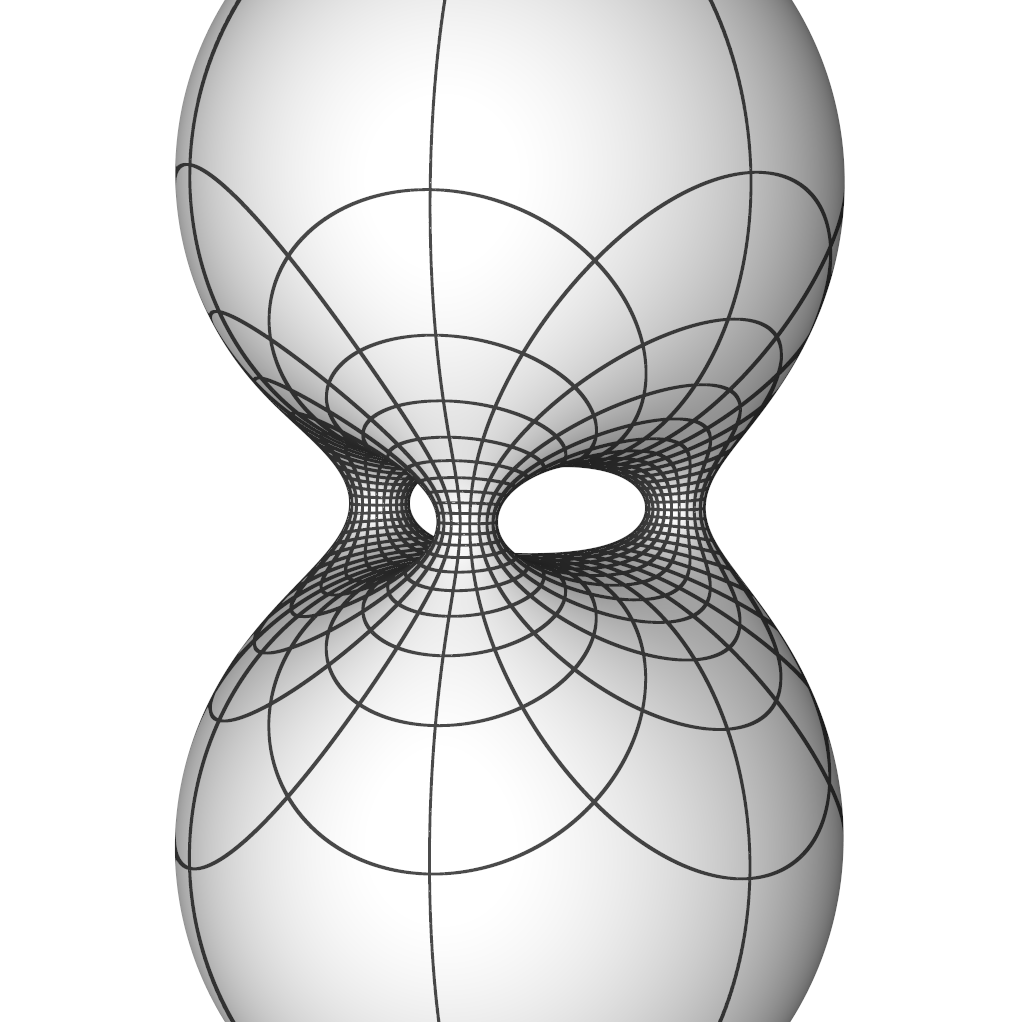}
\includegraphics[width= 0.4\textwidth]{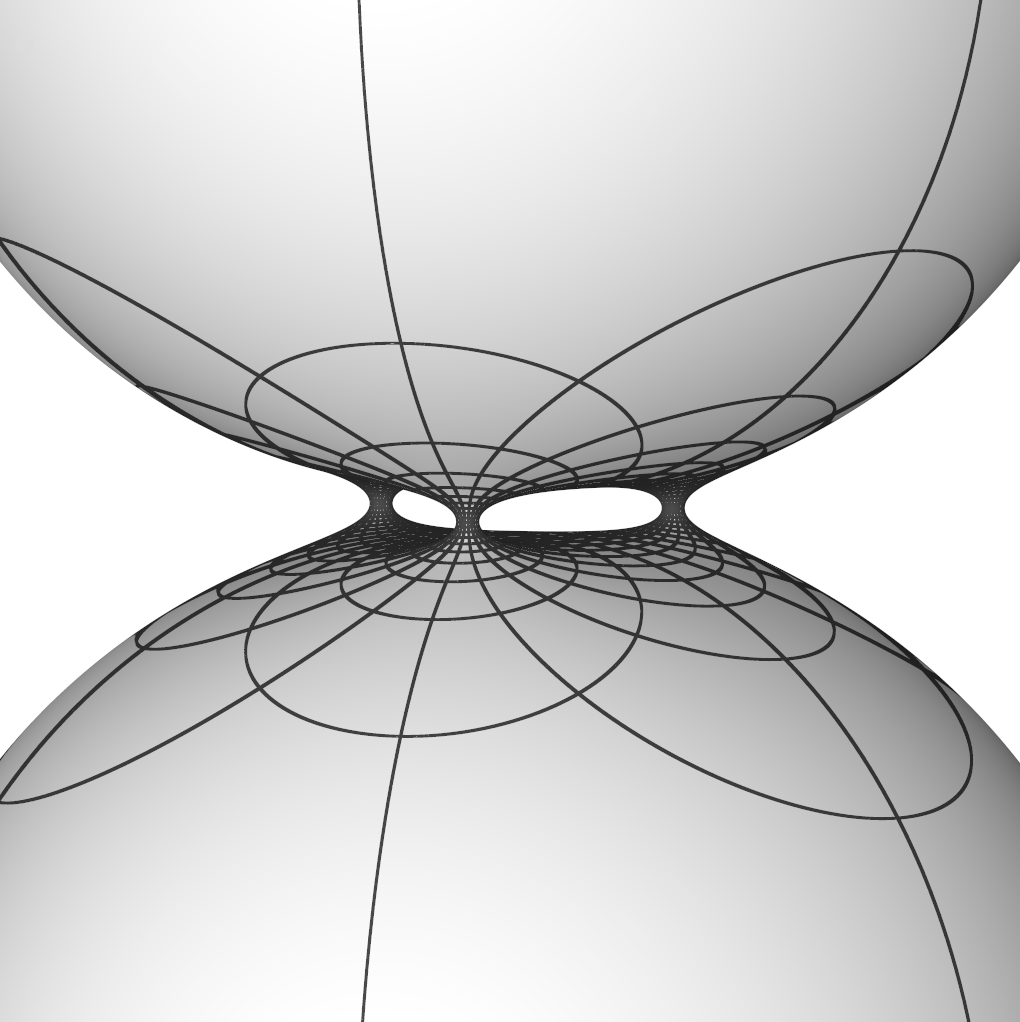}

\caption{ 
Family of embedded constant mean curvature surfaces of genus two in the 3-sphere starting at the Lawson surface $\xi_{2,1}$  and converging to a double cover of an equatorial 2-sphere as the conformal type degenerates \cite{HelSch}. The Willmore energy $\mathcal W$ initially  increases to a maximum above $8\pi$ and then decreases to the limiting value $8\pi$ (see Figure~\ref{fig:graph}).}
\end{figure}

\section{Higher genus outlook}
Examples of compact constrained Willmore surfaces of genus $g\geq 2$ are scarce: the $\Z_{k+1}\times \Z_{l+1}$ symmetric compact minimal surfaces  $\xi_{k,l}$ in the 3-sphere of genus  $g=kl$ found by Lawson \cite{Lawson}, the minimal surfaces in the 3-sphere with Platonic symmetries constructed by Karcher, Pinkall and Sterling \cite{KarcherPinkallSterling}, and the constant mean curvature surfaces of any genus constructed via gluing methods by Kapouleas \cite{Kap}, \cite{Kap2}.
Among Lawson's minimal surfaces Kusner \cite{Kusner}  verified that the family of genus $g$ surfaces $\xi_{g,1}$ have the lowest Willmore energy with $\mathcal W(\xi_{g,1}) < 8\pi$ starting at the Clifford torus $\xi_{1,1}$ with Willmore energy $\mathcal{W}(\xi_{1,1})=2\pi^2$.  The surfaces of Platonic symmetries have larger energies. 
Numerical experiments, using an energy decreasing flow \cite{Kusneretal}, suggest that Lawson's minimal surfaces $\xi_{g,1}$ are the minimizers of $\mathcal W$ over compact surfaces of fixed genus $g$. 

\begin{figure}
\centering
\includegraphics[width=0.4\textwidth]{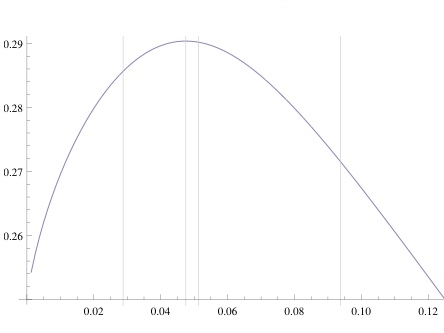}
\includegraphics[width=0.4\textwidth]{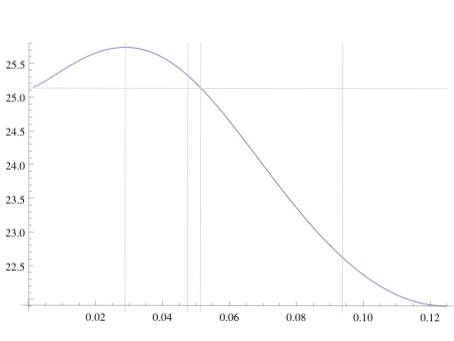}
\caption{Graphs of the mean curvature  $H$ (left) and Willmore energy $\mathcal W$  (right) along the  family of genus two constant mean curvature surfaces (Figure~\ref{fig:highergenus}) deforming Lawson's minimal surface $\xi_{g,1}$. The horizontal axis measures the rectangular conformal type, starting at the square structure for $\xi_{g,1}$ on the right and degenerating  to the twice covered equatorial 2-sphere at the origin. 
}
\label{fig:graph}
\end{figure}

In analogy to the genus one case, where the constrained minimizers in rectangular conformal classes are conjectured to have constant mean curvature (see the 2-lobe Conjecture~\ref{conj:2-lobe}), one could expect such behavior for specific conformal types also for higher genus surfaces. Recently progress has been made in our understanding of the integrable systems approach to higher genus constant mean curvature surfaces \cite{Heller}, \cite{HHS}. Starting from Lawson's minimal surface $\xi_{g,1}$ one can deform $\xi_{g,1}$ by changing the value of the (constant) mean curvature to obtain an (experimentally constructed) 1-parameter family of embedded constant mean curvature surfaces (see Figure~\ref{fig:highergenus}) having  the symmetries of Lawson's $\xi_{g,1}$ minimal surface. The Willmore energy profile  \cite{HHS2} over this family is shown in Figure~\ref{fig:graph}. 
The conformal types of these surfaces are ``rectangular'': the quotient $\P^1$  of the surface under the cyclic $g+1$ fold symmetry has four branch points 
arranged in a rectangle. In particular, the results of \cite{KuwertLi} apply and the existence of a constrained minimizer in those conformal classes for which  $\mathcal W< 8\pi$ is guaranteed.  We also conducted experiments  \cite{HPSW} using the conformal Willmore  flow  \cite{CPS} corroborating the results of  \cite{Kusneretal} and supporting 
\begin{Con}[Higher genus conjecture]
The constrained minimizer of the Willmore energy for a surface of genus $g$ whose conformal type is ``rectangular'' is given by the constant mean curvature analogs of Lawson's minimal surfaces $\xi_{g,1}$. 
\end{Con}
At this juncture one could start speculating  about  the corresponding constrained Willmore Lawson and stability conjectures in higher genus.
It is known \cite{CR} that the space of minimal surfaces in the 3-sphere of fixed genus is compact. A similar result has recently been announced by Meeks and Tinaglia \cite{MT} for strongly
Alexandrov embedded surfaces (these surfaces extend to an immersion of an
embedded compact 3-manifold) of positive constant mean curvature.
 Thus, if one could prove that such surfaces are isolated the correct generalization of Lawson's conjecture (and its constant mean curvature analogue) would be: there are only finitely many (strongly Alexandrov) embedded constant mean curvature surfaces of genus $g$ and fixed mean curvature. They presumably have different conformal types and that would make them the natural candidates for constrained minimizers in those conformal classes. At present time such investigations would be based on very scant evidence due to the lack of generic examples in higher genus. The importance of examples to develop a mathematical  theory has been known throughout its history. In our context this has been well expressed in the opening sentence of Lawson's paper \cite{Lawson}, which has served as a starting point for many of the investigations discussed in this note:  ``It is valuable when dealing with a non-linear theory, such as the study of minimal surfaces, to have available a large collection of examples for reference and insight''.

\end{document}